\documentclass[11pt]{amsart}
\usepackage{calc,amssymb,amsthm,amsmath,fullpage    }
\RequirePackage[dvipsnames,usenames]{xcolor}
\usepackage{hyperref}
\hypersetup{
bookmarks,
bookmarksdepth=3,
bookmarksopen,
bookmarksnumbered,
pdfstartview=FitH,
colorlinks,backref,hyperindex,
linkcolor=Sepia,
anchorcolor=BurntOrange,
citecolor=MidnightBlue,
citecolor=OliveGreen,
filecolor=BlueViolet,
menucolor=Yellow,
urlcolor=OliveGreen
}
\usepackage{alltt}
\usepackage{multicol}
\usepackage{xspace}
\usepackage{rotating}
\usepackage{mabliautoref}
\usepackage{colonequals}
\input{kmacros3.sty}
\usepackage{stmaryrd}

\usepackage{verbatim}
\usepackage{enumerate}
\begin{document}
\title{{Divisor} package for \emph{Macaulay2}}
\author{Karl Schwede}
\date{\today}
\address{Department of Mathematics, University of Utah, 155 S 1400 E Room 233, Salt Lake City, UT, 84112}
\email{schwede@math.utah.edu}

\author{Zhaoning Yang}
\date{\today}							
\email{zyy5054@gmail.com}

\begin{abstract}
This note describes a \emph{Macaulay2} package for handling divisors.  Group operations for divisors are included.  There are methods for converting divisors to reflexive or invertible sheaves.  Additionally, there are methods for checking whether divisors are Cartier, $\bQ$-Cartier, simple normal crossings, generate base point free linear systems, or satisfy numerous other conditions.
\end{abstract}



\thanks{The first named author was supported in part by the
  NSF FRG Grant DMS \#1265261/1501115, NSF CAREER Grant DMS \#1252860/1501102, NSF Grant \#1801849 and a Sloan
  Fellowship.}
\thanks{The second named author was supported in part by the NSF CAREER Grant DMS \#1252860/1501102.}
\maketitle

\section{Introduction}

Divisors are fundamental objects of study within algebraic geometry and commutative algebra.  In this package for \emph{Macaulay2} \cite{M2} we provide a wrapper object for studying Weil and Cartier divisors.  We include tools for studying divisors on both affine and projective varieties.

In this package, divisors are stored (roughly) as formal linear combinations of height one prime ideals, with coefficients from $\bZ$, $\bQ$, or $\bR$.  We include group and scaling operations for divisors, as well as various methods for constructing modules $\O_X(D)$ from divisors $D$ (and vice versa).  We also include code for determining whether divisors are linearly or $\bQ$-linearly equivalent, and for checking whether divisors are Cartier or $\bQ$-Cartier (or finding the non-Cartier locus).  Finally, we also include a number of functions for handling reflexive modules, ideals and their powers.

We realize there is a Divisor class defined in a tutorial in the \emph{Macaulay2} help system.  In that implementation, divisors are given as a pair of ideals---an ideal corresponding to the positive part and an ideal corresponding to the negative part. Our approach offers the advantage that it is easier for the user to see the structure of the divisor.  Additionally, certain operations are much faster in our approach.

We warn the user that when a divisor is created, Gr\"obner bases are constructed for each prime ideal defining a component of the divisor.  Hence, the construction phase may be slower than other potential implementations (and in fact slower than our initial implementation).  However, we feel that this choice offers advantages of execution speed for several functions as well as substantial improvements in code readability.

Within the package, it is tacitly assumed that the ambient ring on which we are working is normal.  This includes the projective case, so care should be taken to make sure the graded ring you are working on satisfies Serre's second condition, see for example \cite[Theorem 8.22A]{Hartshorne} or \cite[Proposition 2.2.21]{BrunsHerzog}.  While one can talk about subvarieties of codimension 1 on more general schemes, the correspondence between divisors and reflexive sheaves is much more complicated, so we restrict ourselves to the normal case.  For an introduction to the theory of rank-1-reflexive sheaves on ``nice'' schemes, see \cite{HartshorneGeneralizedDivisorsOnGorensteinSchemes,HartshonreGeneralizedDivisorsAndBiliaison}; and for a more basic introduction see, for instance, \cite[Chapter II, Sections 5--7]{Hartshorne}.

This paper is structured as follows.  We first give a brief introduction to the construction, conversion, and group operation functions in \autoref{sec.Construction}.  We then discuss the methods for converting divisors $D$ to modules $\O_X(D)$ and converting modules back to divisors in \autoref{sec.Modules}.  \autoref{sec.Checks} describes how to determine if divisors satisfy varies properties (for instance {\tt isCartier} or {\tt isSNC}).  We conclude with a section on future plans.

\subsection*{Acknowledgements}

We thank Tommaso de Fernex, David Eisenbud, Daniel Grayson, Anurag Singh, Greg Smith, Mike Stillman, and the referees for useful conversations and comments on the development of this package.  We also thank the referees for numerous useful comments on this paper.

\section{Construction, conversion and group operations for divisors}
\label{sec.Construction}
This package includes a number of ways to construct a divisor (an object of class {\tt WeilDivisor}), illustrated below.
\begin{verbatim}
i1 : needsPackage "Divisor";
i2 : R = QQ[x,y,u,v]/ideal(x*y-u*v);
i3 : D = divisor({2, 3}, {ideal(x,u), ideal(x, v)})
o3 = 3*Div(x, v) + 2*Div(x, u)
o3 : WeilDivisor on R
i4 : E = divisor(x)
o4 = Div(u, x) + Div(v, x)
o4 : WeilDivisor on R
i5 : F = divisor( (ideal(x,u))^2*(ideal(x,v))^3 )
o5 = 3*Div(v, x) + 2*Div(u, x)
o5 : WeilDivisor on R
\end{verbatim}
The output is a formal sum of height one prime ideals.  The first method requires a list of integers and a list of prime ideals.  The third construction method finds a divisor defined by the given ideal in codimension 1.

We have different classes for $\bQ$-divisors and $\bR$-divisors ({\tt QWeilDivisor} and {\tt RWeilDivisor} respectively), these are constructed via the {\tt divisor} function with the {\tt CoeffType =>} option set or by multiplying a {\tt WeilDivisor} by a rational or real number.  See the documentation.

All types of divisors are ancestors of the {\tt HashTable} class.  Internally, they are hash tables where each key is a list of Gr\"obner basis generators for a prime height-one ideal and each associated value is a list, the first entry of which is the coefficient of the prime divisor and the second entry is the prime ideal used to display the divisor (it tries to match how the user entered it for ease of reading).  Besides the keys corresponding prime divisors, there is a key that specifies the ambient ring and another key that points to a {\tt CacheTable}.

One can convert one type of divisor to another more general class, either by multiplication by appropriate coefficients or by calling appropriate functions.
\begin{verbatim}
i2 : R = QQ[x,y,u,v]/ideal(x*y-u*v);
i3 : D = divisor({1, -3}, {ideal(x,u), ideal(y,u)});
o3 : WeilDivisor on R
i4 : 1/1*D
o4 = -3*Div(y, u) + Div(x, u)
o4 : QWeilDivisor on R
i5 : toQWeilDivisor(D)
o5 = Div(x, u) + -3*Div(y, u)
o5 : QWeilDivisor on R
\end{verbatim}
One can convert $\bQ$ or $\bR$-divisors back to Weil divisors as follows.
\begin{verbatim}
i3 : D = divisor( {2/3, -1/2}, {ideal(x,u), ideal(y, v)}, CoeffType=>QQ)
o3 = 2/3*Div(x, u) + -1/2*Div(y, v) of R
o3 : QDiv
i4 : isWDiv(D)
o4 = false
i5 : isWDiv(6*D)
o5 = true
i6 : toWDiv(6*D)
o6 = 4*Div(x, u) + -3*Div(y, v) of R
o6 : WDiv
\end{verbatim}
See the documentation for more examples.  Alternately, the functions {\tt ceiling} and {\tt floor} will convert any $\bQ$ or $\bR$-divisor to a Weil divisor by taking the ceiling or floor of the coefficients respectively.  More generally, one can call the method {\tt applyToCoefficients} to apply any function to the coefficients of a divisor (since divisors are a type of {\tt HashTable}, this is just done via the {\tt applyValues} function).

Divisors form an Abelian group and one can add {\tt WeilDivisor/QWeilDivisor/RWeilDivisor} to each other to obtain new divisors.  Likewise one can scale by integers, rational numbers or real numbers.
\begin{verbatim}
i3 : D = divisor({1, -2}, {ideal(x,u), ideal(x, v)}); E = divisor(u);
o3 : WeilDivisor on R
o4 : WeilDivisor on R
i5 : 3*D+E
o5 = 4*Div(x, u) + -6*Div(x, v) + Div(u, y)
o5 : WeilDivisor on R
i6 : D - (1/2)*E
o6 = -2*Div(x, v) + 1/2*Div(x, u) + -1/2*Div(u, y)
o6 : QWeilDivisor on R
\end{verbatim}

Since divisors are implemented as subclasses of hash tables, these operations are easily executed internally via the {\tt merge} and {\tt applyValues} commands.

\section{Modules, ideals, divisors and applications}
\label{sec.Modules}

It is well known that divisors are so useful because of their connections with invertible and reflexive sheaves.  This package includes many functions for conversion between these types of objects.  For instance, we have the following.
\begin{verbatim}
1 : R = QQ[x,y,z]/ideal(x*y-z^2); needsPackage "Divisor";
i3 : D = divisor(ideal(x, z));
o3 : WeilDivisor on R
i4 : OO(D)
o4 = image {-1} | x z |
           {-1} | z y |
o4 : R-module, submodule of R
i5 : divisor(o4)
o5 = -Div(z, x)
o5 : WeilDivisor on R
i6 : divisor(o4, IsGraded=>true)
o6 = Div(z, x)
o6 : WeilDivisor on R
\end{verbatim}
The function {\tt OO} produces a module $M$ so that $\widetilde{M} \cong \mathcal{O}_X(D)$ (and the gradings of $M$ are set appropriately).  The function {\tt divisor(M)} only produces a divisor $E$ such that $\O_X(E)$ is isomorphic $\widetilde{M}$.  In particular, ${\tt divisor}({\tt OO}(D))$ will only produce a divisor linearly equivalent to $D$.  

The computation of {\tt OO(D)} is done via a straightforward strategy.  If ${\tt D} = \sum_{i = 1}^m a_i P_i$ where $a_i$ are integers and the $P_i$ are primes, then we can compute $\bigotimes P_i^{-a_i}$ (keeping in mind negative exponents mean applying $\Hom_R(\blank, R)$) and compute the reflexification (see the method {\tt reflexify}).  We do several things make this computation faster.  Firstly, we break up the divisor into the positive and negative parts, and handle them separately (applying the {\tt reflexify} method as little as possible).  Then, instead of computing $P_i^{|a_i|}$, which can have many generators, we form an ideal generated by the generators of $P_i$ raised to the $|a_i|$-th powers.  Since this agrees with $P_i^{|a_i|}$ in codimension 1, it will give the correct answer up to reflexification.  We have noticed substantial speed improvements using this technique.

The function {\tt divisor(Module)} works as follows.  First, it embeds the module as an ideal $I \subseteq R$ via the function {\tt embedAsIdeal}.  After we have an ideal $I$, we call {\tt divisor(I)}.  This finds a divisor $D$ such that $\O_X(D)$ is isomorphic to the given ideal $I$ (in a non-graded sense).  The function {\tt divisor(Ideal)} does this by looking at the minimal height 1 primes $Q_i$ of the ideal $I$ and finding the maximum power $n_i$ such that $I \subseteq Q_i^{(n_i)}$ (the symbolic power).  Note that because $Q_i$ has height 1, we know that $Q_i^{(n_i)} = (Q_i^{n_i})^{**}$ where $\blank^{**}$ denotes reflexification/S2-ification of the ideal.  Finding this maximal power is done by a binary search.  Again, for speed, we compute $(Q_i^{n_i})^{**}$ as $(Q_i^{[n_i]})^{**}$.  If the {\tt IsGraded} flag is set to {\tt true}, {\tt divisor(Module)} corrects the degree of the divisor by adding or subtracting the divisor of an element of appropriate degree (you can see this being done in the example above).  Finding the element of appropriate degree is accomplished via the function {\tt findElementOfDegree}, which uses Smith normal form in the multi-degree setting to solve the system of linear diophantine equations and find a monomial of the given multi-degree.

\begin{remark}
A variant of the function {\tt embedAsIdeal} appeared in the \emph{Macaulay2} documentation in the Divisor tutorial, it also appeared in the work of Moty Katzman.  Our version is slightly more robust than those as it tries to embed the module into the ring in several ways, including some random attempts (see the documentation for how to control the number of random attempts).
\end{remark}

Instead of calling {\tt divisor(Module)}, one can call {\tt divisor(Module, Section => f)}.  This function finds the unique effective divisor $D$ corresponding to a global section $f \in M$ of our module.  The function {\tt divisor(Ideal, Section => f)} behaves similarly.  The strategy is the same as above, additionally one tracks the section and adds a divisor corresponding to the section at the end.

It is worth mentioning that the function {\tt canonicalDivisor} simply computes the canonical module via an appropriate $\Ext$ and then calls {\tt divisor(Module)}.  If you wish to construct a canonical divisor on a projective variety, make sure to set the {\tt IsGraded} option to {\tt true}.

\subsection{Pulling back divisors}

 Utilizing the module and divisor correspondence {\tt pullBack} pulls back a divisor along a map $\Spec S \to \Spec R$ induced by a ring map $R \to S$.  The user has a choice of two algorithms built into this function.  The first works for nearly any map, provided that the divisor is Cartier, and it also works for arbitrary divisors in the flat or finite case.  The second, which is the default strategy, only gives accurate answers if the map is flat, or if the map is finite (or if the prime components of the divisor are Cartier).  It can be faster than the first algorithm, especially for divisors with large coefficients.  To use the first algorithm, use is {\tt Strategy => Sheaves}, to use the second, use the {\tt Strategy => Primes}.

Let us briefly describe these two strategies.  The first algorithm pulls back the sheaf $\O(D)$, keeping track of a section appropriately.
The second algorithm extends each prime ideal defining a prime divisor of $D$ to an ideal of $S$, then it calls {\tt divisor(Ideal)} on each such ideal and sums them keeping track of coefficients appropriately.

Consider the following example where we look at pulling back a divisor after blowing up the origin (we only consider one chart of the blowup).
\begin{verbatim}
i2 : R = QQ[x,y];
i3 : S = QQ[a,b];
i4 : f = map(S, R, {a*b, b});
o4 : RingMap S <--- R
i5 : D = divisor(x*y*(x+y)*(x-y))
o5 = Div(x+y) + Div(-x+y) + Div(x) + Div(y)
o5 : WeilDivisor on R
i6 :  pullback(f, D)
o6 = Div(a+1) + Div(a-1) + 4*Div(b) + Div(a)
o6 : WeilDivisor on S
\end{verbatim}
Note one of the components was lost in this pull-back, as it should have been.  The coefficient of the exceptional divisor is also $4$, as it should be.

\subsection{Global sections}

There are only a few built-in functions for dealing with global sections of modules corresponding to divisors in the current version (in the future we hope to add more tools to do this).  Of course, the user may call {\tt basis(0, OO(D))} to get the global sections of a module corresponding to a divisor.  In this section, we describe briefly two functions for handling global properties of divisors.

The function {\tt mapToProjectiveSpace} gets the global sections of $\O(D)$ and then computes the corresponding map to projective space.  This of course assumes the divisor is graded.  In the example below we project $\mathbb{P}^1 \times \mathbb{P}^1$ to one of its terms by calling {\tt mapToProjectiveSpace} along a divisor of one of the rulings.
\begin{verbatim}
i2 : R = QQ[x,y,u,v]/ideal(x*y-u*v);
i3 : D = divisor(ideal(x,u));
o3 : WeilDivisor on R
i4 : mapToProjectiveSpace(D)
o4 = map(R,QQ[YY , YY ],{v, x})
                1    2
o4 : RingMap R <--- QQ[YY , YY ]
                         1    2
\end{verbatim}

Still assuming the divisor is graded, the function {\tt baseLocus} finds a defining ideal for the locus where $\O(D)$ is \emph{not} generated by global sections.  This is done by computing the cokernel of $\O^{\oplus n} \to \O(D)$ where $H^0(X, \O(D))$ has a basis of $n$ distinct global sections and the map is the obvious one.  In the following example, we compute the base locus of a point on an elliptic curve, and also two times a point on an elliptic curve (which is degree 2 and hence base point free).
\begin{verbatim}
i2 : R = QQ[x,y,z]/ideal(y^2*z-x*(x+z)*(x-z));
i3 : D = divisor( ideal(x,y) );
o3 : WeilDivisor on R
i4 : baseLocus(D)
o4 = ideal (y, x)
o4 : Ideal of R
i5 : baseLocus(2*D)
o5 = ideal 1
o5 : Ideal of R
\end{verbatim}

\section{Checking properties of divisors}
\label{sec.Checks}

The package {\tt Divisor} can check divisors for several properties.  First, we describe the method {\tt isCartier}.
\begin{verbatim}
i2 : R = QQ[x,y,z]/ideal(x^2-y*z);
i3 : D = divisor(ideal(x,y));
i4 : isCartier(D)
o4 = false
i5 : nonCartierLocus(D)
o5 = ideal (z, y, x)
o5 : Ideal of R
i6 : isCartier(2*D)
o6 = true
i7 : isCartier(D, IsGraded => true)
o7 = true
\end{verbatim}
The algorithm behind this function is as follows.  We compute $\O_X(-D) \cdot \O_X(D)$ and check whether it is equal to $\O_X$.  In general, $\O_X(-D) \cdot \O_X(D)$ always defines an ideal defining the non-Cartier locus of $D$, hence the command {\tt nonCartierLocus}.  If the option {\tt IsGraded => true}, then the relevant functions saturate the ideals with respect to the irrelevant ideal.

We also briefly describe the method {\tt isQCartier}.
\begin{verbatim}
i8 : isQCartier(5, D)
o8 = 2
\end{verbatim}
This checks whether any multiples $n \cdot D$ of a Weil divisor or $\bQ$-divisor $D$ are Cartier for any integer $n$ less than or equal to the first argument (in this case $n \leq 5$), it may actually search a little higher than the first argument in the $\bQ$-Cartier case due to rounding issues.  If it finds that $nD$ is Cartier, it returns the integer $n$.  If it doesn't find any Cartier divisors, it returns {\tt 0}.

Some other useful functions are {\tt isPrincipal} and {\tt isLinearEquivalent}.  Checking whether a divisor is principal just comes down to checking whether $\O_X(D)$ is a free module and checking whether $D \sim E$ just boils down to checking whether $D-E$ is principal.  In the graded case, we can do this via \emph{Macaulay2} using the {\tt prune} and {\tt isFreeModule} commands. Unfortunately, we do not know an algorithm for deciding if a non-graded module is free (although we still try to prune the module and more).  Therefore {\tt isPrincipal} and {\tt isLinearEquivalent} can give a false negative for non-graded divisors (the function warns you if this might be the case).  Likewise, the option {\tt IsGraded} can be applied within {\tt isLinearEquivlavent}, which checks that $\O_X(D-E)$ is principal of degree zero.

We can also check whether a divisor $D$ has simple normal crossings by calling {\tt isSNC}.  This first checks that the ambient space of $D$ is regular, then it checks that each prime divisor of $D$  defines a regular scheme, finally it checks that every intersection of of prime divisors of $D$ also defines a regular scheme of the appropriate dimension.

\section{Future plans}
\label{sec.Plans}

There are a number of ways that this package should be expanded.  One of the most important things to be done is to further develop the global methods related to divisors.  We have recently added the ability to check whether a divisor is very ample via the {\tt isVeryAmple} function, which uses the {\tt RationalMaps} package.  However, there is much more to be done.  Some basic intersection theory between divisors and smooth curves would be natural to include.

While the latest version of the package stores the outputs of some functions in the cache, this can still be improved.  For example, there are likely ways to take advantage of knowing that a given divisor is Cartier or $\bQ$-Cartier.

\bibliographystyle{skalpha}
\bibliography{MainBib}

\end{document}